\documentclass[1p]{elsarticle}
\usepackage{amsthm,amsmath,amssymb,amsfonts}
\usepackage{xcolor}

\newtheorem{theorem}{Theorem}
\newtheorem{conjecture}[theorem]{Conjecture}

\newtheorem{problem}[theorem]{Problem}
\newtheorem{claim}{Claim}

\newproof{pf}{Proof}

\begin{document}

\title{The irregularity strength of dense graphs -- on asymptotically optimal solutions of problems of Faudree, Jacobson, Kinch and Lehel}

\author[agh]{Jakub Przyby{\l}o%\corref{cor1}
\fnref{MNiSW}
\fnref{DecInt}}
\ead{jakubprz@agh.edu.pl, phone: 048-12-617-46-38} %,  fax: 048-12-617-31-65}

\fntext[MNiSW]{Research supported by AGH University of Krakow 16.16.420.054.}
%This work was partially supported by the Faculty of Applied Mathematics AGH UST statutory tasks within subsidy of Ministry of Science and Higher Education.}
\fntext[DecInt]{Declarations of interest: none.}

\address[agh]{AGH University of Krakow, Faculty of Applied Mathematics, al. A. Mickiewicza 30, 30-059 Krakow, Poland}

\begin{abstract}
The irregularity strength of a graph $G$, $s(G)$, is the least $k$ such that there exists a $\{1,2,\ldots,k\}$-weighting of the edges of $G$ attributing distinct weighted degrees to all vertices, or equivalently the least $k$ 
enabling obtaining a multigraph with nonrecurring degrees by blowing each edge $e$ of $G$ to at most $k$ copies of $e$. In 1991 Faudree, Jacobson, Kinch and Lehel asked for the optimal lower bound for the minimum degree of a graph $G$ of order $n$ which implies that $s(G)\leq 3$. More generally, they also posed a similar 
question regarding the upper bound $s(G)\leq K$ for any given constant $K$. We provide asymptotically tight solutions of these problems by proving that such optimal lower bound is of order $\frac{1}{K-1}n$ for every fixed integer $K\geq 3$.
\end{abstract}


\maketitle

\section{The irregularity strength of graphs\label{Subsection_Origins}}

Suppose we call a graph (respectively, multigraph) \emph{irregular} if all its vertices have pairwise distinct degrees. Though this seems the most natural antonym of a \emph{regular} graph, it is the most basic fact that
every nontrivial graph is not irregular in this sense. This observation motivated research on several alternative 
notions discussed by Chartrand, Erd\H{o}s and Oellermann in the paper~\cite{ChartrandErdosOellermann} entitled \emph{How to define an irregular graph?}. Other notions of this type were also discussed e.g. in~\cite{IrregularGraphs3,IrregularGraphs2}, but none of these seemed to outshine the remaining ones.
In~\cite{Chartrand} Chartrand et al. thus turned towards measuring the `irregularity' of a graph, exploiting in particular the fact that there are irregular multigraphs of any order larger than $2$. 
Given any (simple) graph $G=(V,E)$ which does not contain an isolated edge nor more than one isolated vertex, they defined the \emph{irregularity strength} of $G$, $s(G)$, as the least $k$ so that one may obtain an irregular multigraph of $G$ by blowing each edge $e$ of $G$ into at most $k$ copies of $e$. We shall call such graphs \emph{good}.
We set $s(G)=+\infty$ for the remaining ones. This definition can be conveniently rephrased by means of $k$-weightings, i.e. mappings $\omega:E\to[k]$, where $[k]=\{1,2,\ldots,k\}$. 
Then $s(G)$ is simply the least $k$ admitting such a $k$-weighting for which each vertex $v$ of $G$ is attributed a distinct \emph{weighted degree} $\sigma_\omega(v):=\sum_{u\in N(v)}\omega(uv)$, where $N(v)=N_G(v)$ is the set of neighbours of $v$ in $G$; we shall also call $\sigma_\omega(v)$ the \emph{sum of} or \emph{at} $v$.
It is known that $S(G)\leq n-1$ for any good graph of order $n$ which is not the triangle, see~\cite{Aigner,Nierhoff}. This upper bound is sharp e.g. for the family of stars. Much better bounds can however be achieved for graphs with larger minimum degree $\delta$. 
In particular, in~\cite{KalKarPf} Kalkowski, Karo\'nski and Pfender proved that for every good graph $G$ with $\delta>0$, 
\begin{equation}\label{6ndeltaBound}
s(G)\leq 6\left\lceil \frac{n}{\delta}\right\rceil.
\end{equation}
This upper bound is tight up to a multiplicative factor, as exemplified by $d$-regular graphs, for which
a straightforward counting argument, cf.~\cite{Chartrand}, implies that
\begin{equation}\label{LoweBoundRegular}
s(G)\geq \frac{n+d-1}{d} = \frac{n}{d}+1-\frac{1}{d}.
\end{equation}
In fact this lower bound was one of the main motivations behind the following key open conjecture in this field, posed
 in 1987 by Faudree and Lehel. (Actually, it was first asked as a  question by Jacobson, as Lehel mentions in~\cite{Lehel}.)
\begin{conjecture}[Faudree and Lehel~\cite{Faudree}]\label{FaudreeAndLehelConjecture}
There exists an absolute constant $C$ such that for every $d$-regular graph $G$ with $n$ vertices and $d\geq 2$,
$$s(G)\leq\frac{n}{d}+C.$$
\end{conjecture}
As a matter of fact this was the problem that brought considerable attention and
 “energised the study of the irregularity strength”, as Cuckler and Lazebnik  affirm in~\cite{Lazebnik}.
It resulted in an extensive list of papers devoted to this graph invariant, see in particular~\cite{Aigner,Amar,Amar_Togni,Bohman_Kravitz,Lazebnik,Dinitz,Ebert,Ebert2,Faudree2,Faudree,Ferrara,Frieze,Gyarfas,Jendrol_Tkac,KalKarPf,Lehel,MajerskiPrzybylo2,Nierhoff,PrzybyloIrregRandomGr,AsymptoticIrregStrReg,Przybylo,irreg_str2,PrzybyloWeiLong,PrzybyloWeiShort,Togni}.
This parameter can moreover be regarded a cornerstone of entire rich branch of graph theory, as it triggered research devoted to vastness of related concepts, see 
\cite{Louigi30,Louigi2,Louigi,BarGrNiw,LocalIrreg_1,BensmailMerkerThomassen,BonamyPrzybylo,
Hatami,Joret,KalKarPf_123,123KLT,Keusch2,Keusch1,
Przybylo_asym_optim,LocalIrreg_2,Przybylo_CN_1,1234Reg123,
Przybylo_CN_2,ThoWuZha,Vuckovic_3-multisets,WongZhu23Choos,WongZhuChoos} 
for examples of papers regarding such problems, many of which can also be found in Gallian's survey~\cite{Gallian_survey}.
It is believed that a natural generalisation of Conjecture~\ref{FaudreeAndLehelConjecture} towards all graphs, with $d$ replaced by the minimum degree $\delta$, should hold. 
Nevertheless, it was a long standing open problem if there was at least is a linear in $n/d$ ($n/\delta$, resp.), upper bound for $s(G)$.
This was partly confirmed by Frieze, Gould, Karo\'nski and Pfender~\cite{Frieze} and by Cuckler and Lazebnik~\cite{Lazebnik}, to be later finally settled in~\cite{Przybylo,irreg_str2}, and then improved to~\eqref{6ndeltaBound} in~\cite{KalKarPf}.
Quite recently Przyby{\l}o showed that 
$s(G) \leq \frac{n}{d}(1+ \frac{1}{\ln^{\varepsilon/19}n})=\frac{n}{d}(1+o(1))$ 
for $d$-regular graphs with 
$d\in [\ln^{1+\varepsilon} n, n/\ln^{\varepsilon}n]$.
This bound was then improved and extended towards all values of $d$ by Przyby{\l}o and Wei~\cite{PrzybyloWeiLong,PrzybyloWeiShort}, and generalised to all good graphs with minimum degree $\delta\geq 1$ in~\cite{PrzybyloWeiLong}, 
which implies in particular that $s(G)\leq \frac{n}{\delta}(1+O(\delta^{-\varepsilon}))$ for any fixed $\varepsilon \in (0,1/4)$ and moreover includes the following result.
\begin{theorem}[\cite{PrzybyloWeiLong}]\label{thm:main1-long_paper}
For every $\beta > 0.8$, there is an absolute constant $C$ such that for each graph $G$ with $n$ vertices and minimum degree 
$\delta\geq n^{\beta}$,  
\[s(G) \leq \frac{n}{\delta} + C.\]
\end{theorem}

\section{Dense graphs}

In 1991 Lehel published an undoubtedly essential for development of this field  survey~\cite{Lehel} entitled
\emph{Facts and quests on degree irregular assignments}, including an array of key problems related with 
the irregularity strength of graphs. Tellingly, many of these still remain open.
In particular the mentioned Conjecture~\ref{FaudreeAndLehelConjecture}, 
which was included in the chapter \emph{Regular and Dense Graphs} of~\cite{Lehel}.
The second part of this chapter contained several problems concerning dense graphs, or more precisely -- graphs with large minimum degree, which were supposed to have irregularity strength bounded above by a constant. 
The base question here concerns the condition implying that $s(G)\leq 3$. 
\begin{problem}[Problem 29,~\cite{Lehel}]\label{Problem29}
Find the largest integer $t=t(n)$ such that $s(G)\leq 3$ holds, if $n$ is sufficiently large, for every graph $G$ of order $n$ with minimum degree $n-t$.
\end{problem}
By~\cite{Chartrand}, $s(K_n)=3$, and hence $3$ is indeed the least upper bound which may be considered in context of such a question. Problem~\ref{Problem29} was in fact posed by Faudree, Jacobson, Kinch and Lehel~\cite{FaudreEtAlDenseIrregular}, who also showed that such non-trivial $t$ exists.
\begin{theorem}[\cite{FaudreEtAlDenseIrregular}]\label{DenseIrregularStrGraphsTh}
Let $G$ be a graph of order $n$ and minimum degree $\delta = n-t$. If $1\leq t \leq \sqrt{\frac{n}{18}}$, then $s(G)\leq 3$.
\end{theorem}
Note that by \eqref{LoweBoundRegular}, $s(G)> 3$ for 
every integer $d<\lfloor n/2\rfloor$. Therefore, we must have $t(n)\leq \lceil n/2\rceil$ for the function $t$ in Problem~\ref{Problem29}.
The same observation yields also the base upper bound (for $t$) 
in the following dual of this problem, posed in~\cite{FaudreEtAlDenseIrregular} and included in Lehel's survey.
\begin{problem}[Problem 31,~\cite{Lehel}]\label{Problem31}
Find the smallest integer $t=t(n)$ such that there exists infinitely many graphs of order $n$ with minimum degree $n-t$ and having irregularity strength at least $4$.
\end{problem}
More generally, Faudree, Jacobson, Kinch and Lehel~\cite{FaudreEtAlDenseIrregular} posed the problem concerning any constant upper bound for the irregularity strength.
\begin{problem}[Problem 26,~\cite{Lehel}]\label{Problem26}
Let $\alpha$ be a real with $0 < \alpha < 1$. Does there exist a constant $K=K(\alpha)$
 such that $s(G)\leq K$ for every graph $G$ of order $n$ with minimum degree $\delta \geq (1-\alpha)n$?
\end{problem}
This issue was first solved in the affirmative in the mentioned paper of Cuckler and Lazebnik~\cite{Lazebnik}, who proved that $s(G)\leq 48n/\delta +6$ for graphs with  minimum degree $\delta\geq 10n^{3/4}\ln^{1/4}n$. 
Their bound was then improved  in~\cite{KalKarPf} to~\eqref{6ndeltaBound}, which implies that 
 $s(G)< 6\frac{n}{\delta} +6\leq  \frac{6}{1-\alpha}+6$ for graphs with  $\delta\geq (1-\alpha)n$.
 However, this bound is roughly $6$ times larger than the lower bound stemming from~\eqref{LoweBoundRegular}.
At the same time, $ \frac{6}{1-\alpha}+6\geq 12$ for every $\alpha\in(0,1)$, so this result cannot imply $K(\alpha)$ smaller than $12$ in Problem~\ref{Problem26}. 
In fact, this problem can be analysed from a slightly different angle. Since $s(G)$ can only take integer values, it asks for the thresholds for $\delta$ (in terms of $n$) above which $s(G)$ must be below given integer values. From this perspective, Problem~\ref{Problem26} is a direct generalisation of Problem~\ref{Problem29}, which regards the first such threshold.  
As our main result, we shall provide asymptotically optimal solutions of these problems for all $K\geq 3$.
\begin{theorem}\label{main_theorem_constant_K}
For every fixed $K\geq 3$, there exists a constant $C$ such that $s(G)\leq K$ for every graph $G$ of order $n\geq 2$ with minimum degree 
\begin{equation}\label{DeltaLowerInMain}
\delta\geq \left(\frac{1}{K-1}+\frac{C}{\left(\frac{n}{\ln n}\right)^\frac{1}{6}}\right)n.
\end{equation}
\end{theorem}
Note that given any integer $K\geq 3$ it is easy to construct a $d$-regular graph $G$ with $n=d(K-1)+2$ vertices (there are many such graphs, for arbitrarily large $n$). By~\eqref{LoweBoundRegular} however, $s(G)> K$, while
\begin{equation}\label{LowerBoundAboveK}
\delta(G)=d=\frac{n-2}{K-1} = \left(\frac{1}{K-1}-\frac{2}{(K-1)n}\right)n = \left(\frac{1}{K-1}-o(1)\right)n.
\end{equation}
Therefore, the lower bound for $\delta$ in Theorem~\ref{main_theorem_constant_K} is indeed asymptotically optimal.

Finally, let us mention that Theorem~\ref{thm:main1-long_paper} implies e.g. that for any fixed $0.8<\beta<1$ there is a (large) constant $C$ such that $s(G)\leq K$ if $\delta\geq \frac{n}{K-C}+n^\beta = (\frac{1}{K-C}+o(1))n$. This bound is however significantly worst than the one in  Theorem~\ref{main_theorem_constant_K}, by a multiplicative constant factor, for every fixed $K$, and moreover applies only to $K>C$.

\section{Tools and Notation}

Given a graph $G=(V,E)$, $v\in V$, $U\subseteq V$ and $F\subseteq E$, by $d_U(v)$ we denote the number of $u\in U$ such that $uv\in E$, and by $d_F(v)$ we mean the number of edges $uv\in F$.
We also denote by $G[U]$ the subgraph induced by $U$ in $G$.

We shall be using the following one-sided Chernoff Bounds, as well as its two-sided variant without restriction on the value of deviation $t$.

\begin{theorem}[Chernoff Bound, cf., e.g.,~\cite{MolloyReed}]\label{ChernoffBoundTh1}
Let $X_1, \dots, X_n$ be i.i.d.~random variables such that $\Pr(X_i=1) = p$ and $\Pr(X_i=0)=1-p$. Then for any
 $0\leq t\leq np$,
\begin{equation}\label{DoubleChernoffBound}
\mathbf{Pr} \left(\sum_{i=1}^n X_i >np+t\right)<e^{-\frac{t^2}{3np}}~~{and}~~\mathbf{Pr}\left(\sum_{i=1}^n X_i <np-t\right)<e^{-\frac{t^2}{2np}}\leq e^{-\frac{t^2}{3np}}.
\end{equation}
\end{theorem}

\begin{theorem}[Chernoff Bound, cf., e.g., \cite{ASbook}, Appendix A]\label{bound:chernoff2}
Let $X_1, \dots, X_n$ be i.i.d.~random variables such that $\Pr(X_i=1) = p$ and $\Pr(X_i=0)=1-p$. Then for any $t \geq 0$,
\begin{align*}
\mathbf{Pr}\left(\left|\sum_{i=1}^n X_i - np\right| > t\right) \leq 2e^{-\frac{t^2}{3\max\{np, t\}}}.
\end{align*} 
\end{theorem}

Note the first inequality in~\eqref{DoubleChernoffBound}
remains valid if $\mathbf{Pr}(X_i=1)\leq p$ for each $i$, 
while the second one -- when $\mathbf{Pr}(X_i=1)\geq p$ for every $i$.

\section{Proof of Theorem~\ref{main_theorem_constant_K}}

\subsection{General idea}

Assume $K\geq 3$ is a fixed constant and $G=(V,E)$ is a graph of order $n$ and with minimum degree $\delta$. We shall prove that if $n$ is large enough, i.e. larger than some constant $n_K$, and 
$\delta\geq (\frac{1}{K-1}+C_1\frac{\ln^{1/6}n}{n^{1/6}})n$ 
 for some fixed, appropriately chosen constant $C_1$, then $s(G)\leq K$. Note this shall imply Theorem~\ref{main_theorem_constant_K}.
 
Below, by $d(v)$ we always mean $d_G(v)$.

In order to prove that $s(G)\leq K$ we shall be using only integers $1$, $2$, $3$ and $K$ to weight the edges.
Our argument can be divided to three phases, which we may roughly outline as follows.
 In the firs stage  we shall randomly partition $V$ to a small set $S$ and a big set $B$. With every vertex we shall also associate a random variable $X_v$ uniformly distributed over $[0,1]$ interval.
Most of the edges in $G$ shall be then initially weighted according to a simple rule: assign weight $\omega_0(uv)=K$ if $X_u+X_v\geq 1$, and assign $\omega_0(uv)=1$ otherwise. As each vertex $v$ shall thus have tendency towards roughly $X_v$ fraction of its incident edges weighted $K$ and $X_v$ shall be independent for all $v\in V$, this shall admit to distribute the initial sums (weighted degrees) $\sigma_0(v)$ for all $v\in V$ relatively sparsely over $\mathbb{N}$. Note that if we denote by $L$ the set of edges weighted $K$, then we shall in particular have $\sigma_0(v)=Kd_L(v)+d_{E\smallsetminus L}(v)$.

Not all edges shall be initially weighted according to the simple rule above. Specifically, the edges in $S$ shall all be weighted $1$ and some edges between $B$ and $S$ shall be (at random) placed in a set $M$, whose members shall all be weighted $1$ as well. Within Claim~\ref{MainRandomClaim}, handling the randomised part of our construction, we shall however in a way control and limit the number of such edges. 
 We shall in particular guarantee that each vertex $v$ in $B$ is incident with a large enough number of such edges, so that within the second stage of our construction we are able to make sure the sum at $v$ is distinct from all other sums (in $B$) by altering the weights of some edges in $M$ from $1$ to $2$.
 
The mentioned Claim~\ref{MainRandomClaim} shall in fact enable us to show that in any interval $I_j$, to be specified later, the initial sums are distributed sparsely enough so that we are able to avoid some integers as sums in $B$ while adjusting these in the mentioned stage 2. This shall facilitate the later adjustments of the sums in $S$, which shall get different from those in $B$, as we shall have more than necessary available `candidates' for their values. The adjustment of sums in $S$ shall be based on changing weights of edges in $S$ from $1$ to $2$ or $3$, and its successful feasibility shall stem mainly from two additional facts. Firstly, sums in $S$ shall be very sparsely distributed, due to (random) stage 1 (and the fact that $S$ shall be `small'). Secondly, each vertex in $S$ shall be incident with large enough number of edges in $S$ (also due to Claim~\ref{MainRandomClaim}), which shall allow a special adjusting algorithm to successively operate  
in stage 3.

\subsection{Initial random weighting}

Rather than explicitly setting the lower bound for $n$, for our fixed $K\geq 3$, we shall be assuming further on that $n$ is large enough so that all the statements and inequalities below hold. Let us set
$$t:=\left(\frac{n}{\ln n}\right)^{\frac13}.$$

Assume  
\begin{equation}\label{deltaBound}
\delta\geq \frac{n}{K-1} + \frac{29Kn}{\sqrt{t}}.
\end{equation}

Set $$i_j:=\delta+j\frac{n}{t}$$ for $j=0,1,\ldots,Kt$, and

$$I_j:=[i_{j-1},i_j)$$ for $j=1,2,\ldots,Kt$.
These intervals shall be required to contain a limited number of initial sums of vertices. For technical reasons we shall however have to apply a two-stage approximation of the densities of initial sums.
Namely, we shall on one hand be able to control the densities of some values $S_v$ related with vertices $v\in V$, which shall be determined by certain random variables $X_v,Y_v$ associated with $v\in V$.
On the other hand, we shall assure that the actual initial sums, resulting from Claim~\ref{MainRandomClaim} below, are at distance at most $3Kn/t^{3/2}$ from the mentioned values $S_v$. Wherefore, we shall need to bound the number of $S_v$ in a slightly wider auxiliary intervals, defined as
$$I'_j:=\left[i_{j-1}-\frac{3Kn}{t\sqrt{t}},i_j+\frac{3Kn}{t\sqrt{t}}\right)$$  
for $j=1,2,\ldots,Kt$.

\begin{claim}\label{MainRandomClaim}
We may associate with every vertex $v\in V$ a real number $x_v\in[0,1]$ and an integer $y_v\in\{0,1\}$
and choose a subset $M$ of the set $\{uv\in E: y_u+y_v = 1\}$ 
and a set  $L\subseteq E\smallsetminus M$ 
such that if we denote
\begin{eqnarray}
B&:=&\left\{v\in V: y_v=0\right\}, \nonumber\\
S&:=&\left\{v\in V: y_v=1\right\}, \nonumber\\
\sigma_v&:=& Kd_L(v)+d_{E\smallsetminus L}(v), \nonumber\\
S_v&:=&d(v)+(K-1)x_vd(v)\left(1-\frac{2Ky_v}{\sqrt{t}}+\frac{(2K-1)(2y_v-1)}{t}\right), \nonumber
\end{eqnarray}
for every $v\in V$,
then for each $v\in V$ and every $j=1,2,\ldots,Kt$:
\begin{enumerate}
\item[$A_{v}$: ] $|\sigma_v-S_v|\leq \frac{3Kn}{t\sqrt{t}}$;
\item[$A_{B,j}$: ] $\left|\left\{v\in B: S_v\in I'_j\right\}\right| \leq \left(1-\frac{17K^2}{\sqrt{t}}\right)\frac{n}{t}$;
\item[$A_{S,j}$: ] $\left|\left\{v\in S: S_v\in I'_j\right\}\right| \leq \left(1-\frac{1}{\sqrt{t}}\right)\frac{n}{t\sqrt{t}}$; 
\item[$A_{B,v}$: ] if $v\in B$, then $d_M(v)\geq 2\frac{n}{t}$; 
\item[$A_{S,v}$: ]  if $v\in S$, then $d_M(v)\leq \frac{2Kd(v)}{\sqrt{t}}$; 
\item[$A'_{S,v}$: ]  if $v\in S$, then $\left|d_S(v)-\frac{d(v)}{\sqrt{t}}\right|\leq \sqrt{n}$. 
\end{enumerate}
\end{claim}

\begin{pf}
Associate with vertices $v\in V$ and edges $e\in E$ independent:
\begin{itemize}
\item random variables $X_v\sim U[0,1]$;
\item binary random variables $Y_v$ taking $1$ with probability $1/\sqrt{t}$ and $0$ otherwise;
\item random variables $Z_e$ taking $1$ with probability $(2K-1)/\sqrt{t}$ 
and $0$ otherwise.
\end{itemize}
The values of $X_v$ shall be used to choose a subset $L$ of $E$ (whose element shall be weighted $K$), the variables $Y_v$ shall indicate a partition of $V$ to $S$ and $B$, while the variables $Z_e$ shall be accountable for marking out the elements of $M$ (among edges joining $S$ and $B$). More specifically,  
we set:
\begin{eqnarray}
B&=&\left\{v\in V: Y_v=0\right\}, \nonumber\\
S&=&\left\{v\in V: Y_v=1\right\}, \nonumber\\
M&=&\left\{uv\in E: Y_u+Y_v=1\wedge Z_{uv}=1\right\}, \nonumber\\
L&=&\left\{uv\in E: X_u+X_v\geq 1 \wedge Y_u+Y_v\leq 1\right\}\smallsetminus M. \nonumber
\end{eqnarray}
We shall show that with positive probability there exist values $x_v$, $y_v$, $z_e$ for all $X_v$, $Y_v$, $Z_e$, $v\in V$, $e\in E$, respectively, such that all our requirements are fulfilled. Let us remark that within our proof we shall be abusing slightly the fixed notation of $S_v$, which shall be considered as a random variable, dependent not on specific $x_v$ and $y_v$, but on the random variables $X_v$ and $Y_v$, respectively (unless their values are fixed as given $x_v$ and $y_v$).

For each of the specified $6$ types of events, $A_v$ -- $A'_{S,v}$, we shall show that its non-occurrence is highly unlikely, i.e. has probability at most $n^{-2}$. 
Let us fix arbitrary $v\in V$. 
Note that the probability that an edge $uv\in E$ is included in $M$ provided that $v\in B$
equals exactly $\frac{1}{\sqrt{t}} \cdot\frac{2K-1}{\sqrt{t}}=\frac{2K-1}{t}$, and thus, by~\eqref{deltaBound},
$$\mathbf{E}\left(d_M(v)~|~ v\in B\right) = \frac{(2K-1)d(v)}{t} 
> \frac{(2K-1)\frac{n}{K-1}}{t} = \frac{2n}{t}+\frac{n}{(K-1)t}.$$ Therefore, by the Chernoff Bound, i.e. Theorem~\ref{ChernoffBoundTh1} (for $n$ large enough),
\begin{equation}\label{ProbBound1}
\mathbf{Pr}\left(\overline{A_{B,v}}\right)\leq \mathbf{Pr}\left(d_M(v)< 2\frac{n}{t}~|~v\in B \right)
\leq e^{-\frac{\left(\frac{n}{(K-1)t}\right)^2}{2\left(\frac{2n}{t}+\frac{n}{(K-1)t}\right)}}<n^{-2}.
\end{equation}

Analogously, the probability that an edge $uv\in E$ is included in $M$ provided that $v\in S$ equals  $(1-\frac{1}{\sqrt{t}}) \cdot \frac{2K-1}{\sqrt{t}} < \frac{2K-1}{\sqrt{t}}$, and thus $\mathbf{E}(d_M(v)~|~ v\in S) < \frac{(2K-1)d(v)}{\sqrt{t}}$. Therefore, by the Chernoff Bound,
\begin{equation}\label{ProbBound2}
\mathbf{Pr}\left(\overline{A_{S,v}}\right)\leq \mathbf{Pr}\left(d_M(v)> \frac{2Kd(v)}{\sqrt{t}}~|~v\in S \right)
\leq e^{-\frac{\left(\frac{d(v)}{\sqrt{t}}\right)^2}{3\frac{(2K-1)d(v)}{\sqrt{t}}}}<n^{-2}.
\end{equation}

Next, since
$\mathbf{E}(d_S(v)~|~v\in S) = \frac{d(v)}{\sqrt{t}}$, then by the Chernoff Bound,
%Theorem~\ref{bound:chernoff2},
\begin{equation}\label{ProbBound3}
\mathbf{Pr}\left(\overline{A'_{S,v}}\right)\leq \mathbf{Pr}\left(\left|d_S(v)-\frac{d(v)}{\sqrt{t}}\right|>\sqrt{n}~|~v\in S \right)
\leq 2e^{-\frac{n}{3\frac{d(v)}{\sqrt{t}}}}<n^{-2}.
\end{equation}

Further, in order to prove that $\mathbf{Pr}(\overline{A_v})\leq n^{-2}$,
it is sufficient to show that $\mathbf{Pr}(\overline{A_v}~|~X_v=x_v, Y_v=y_v)\leq n^{-2}$ for any possible values $x_v,y_v$. Thus let us fix arbitrary $x_v\in[0,1]$, $y_v\in\{0,1\}$. Note that $\sigma_v = d(v)+(K-1)d_L(v)$. Therefore, it is straightforward to compute (it is convenient to separately consider two cases: $y_v=0$ and $y_v=1$) that: 
\begin{eqnarray}
&&\mathbf{E}\left(\sigma_v~|~X_v=x_v, Y_v=y_v\right) \nonumber\\
&=& d(v)+(K-1) ~\mathbf{E}\left(d_L(v)~|~X_v=x_v, Y_v=y_v\right) \nonumber\\
&=& d(v)+(K-1) d(v) x_v 
\left[1-y_v\frac{1}{\sqrt{t}} -y_v\left(1-\frac{1}{\sqrt{t}}\right)\frac{2K-1}{\sqrt{t}}
-\left(1-y_v\right)\frac{1}{\sqrt{t}}\frac{2K-1}{\sqrt{t}}\right] \nonumber\\
&=& d(v)+(K-1)x_vd(v)\left(1-\frac{2Ky_v}{\sqrt{t}}+\frac{(2K-1)(2y_v-1)}{t}\right) \nonumber\\
&=& S_v. \nonumber
\end{eqnarray}

In fact, by the same calculations, for any fixed values $x_v, y_v$ of $X_v, Y_v$, respectively,  $d_L(v)=(\sigma_v-d(v))/(K-1)$ is a sum
of $d(v)$ independent random variables (related with neighbours $u$ of $v$, and determined by $X_u$, $Y_u$ and $Z_{uv}$) taking value $1$ with probability 
$$x_v\left(1-\frac{2Ky_v}{\sqrt{t}}+\frac{(2K-1)(2y_v-1)}{t}\right)$$ 
and $0$ otherwise. Thus, by the Chernoff Bound in Theorem~\ref{bound:chernoff2},
\begin{eqnarray}
\mathbf{Pr}\left(\overline{A_v}~|~X_v=x_v, Y_v=y_v\right) 
&\leq& 
2e^{-\frac{\left(\frac{3Kn}{(K-1)t\sqrt{t}}\right)^2}{3\max\left\{\mathbf{E}\left(d_L(v)~|~X_v=x_v, Y_v=y_v\right),\frac{3Kn}{(K-1)t\sqrt{t}}\right\}}} \nonumber\\
&\leq& 2e^{-\frac{\left(\frac{3K}{K-1}\right)^2n\ln n}{3d(v)}} \nonumber\\
&\leq& 2e^{-\frac{\left(\frac{3K}{K-1}\right)^2\ln n}{3}}<n^{-2}. \nonumber
\end{eqnarray}
Therefore,
\begin{equation}\label{ProbBound4}
\mathbf{Pr}\left(\overline{A_v}\right) < n^{-2}. 
\end{equation}

For any fixed $j\in\{1,2,\ldots,Kt\}$, set 
$$B_j = \{u\in B: S_u\in I'_j\}.$$
Note that the random variable $S_v$, provided that $Y_v=0$,
 is uniformly distributed over an interval of length 
$(K-1)d(v)(1-\frac{2K-1}{t})\geq (K-1)\delta(1-\frac{2K-1}{t})$.
 Therefore, $|B_j|$ is a sum of $n$ independent random variables taking value $1$ with probability at most 
$$\frac{\left(1-\frac{1}{\sqrt{t}}\right)|I'_j|}{(K-1)\delta \left(1-\frac{2K-1}{t}\right)}
\leq \frac{\left(1-\frac{1}{\sqrt{t}}\right)\left(\frac{n}{t}+\frac{6Kn}{t\sqrt{t}}\right)}{(K-1) \left(\frac{n}{K-1} + \frac{29Kn}{\sqrt{t}}\right) \left(1-\frac{2K-1}{t}\right)}
$$ 
and $0$ otherwise. Hence,
\begin{eqnarray}
\mathbf{E}(|B_j|) &\leq& n \cdot \frac{\left(1-\frac{1}{\sqrt{t}}\right)\left(\frac{n}{t}+\frac{6Kn}{t\sqrt{t}}\right)}{\left(K-1\right) \left(\frac{n}{K-1} + \frac{29Kn}{\sqrt{t}}\right) \left(1-\frac{2K-1}{t}\right)} \nonumber\\ 
&=& \frac{\left(1-\frac{1}{\sqrt{t}}\right)\left(1+\frac{6K}{\sqrt{t}}\right)}{ \left(1 + \frac{29K\left(K-1\right)}{\sqrt{t}}\right) \left(1-\frac{2K-1}{t}\right)} \cdot \frac{n}{t} 
= \frac{1+\frac{6K-1}{\sqrt{t}}+o\left(\frac{1}{\sqrt{t}}\right)}{1+\frac{29K\left(K-1\right)}{\sqrt{t}}+o\left(\frac{1}{\sqrt{t}}\right)} \cdot \frac{n}{t}\nonumber\\
 &<&  \left(1-\frac{17K^2+1}{\sqrt{t}}\right) \cdot \frac{n}{t} 
 = \left(1-\frac{17K^2}{\sqrt{t}}\right) \cdot \frac{n}{t} - \frac{n}{t\sqrt{t}}.\nonumber
\end{eqnarray}
Thus, by the Chernoff Bound, 
\begin{equation}\label{ProbBound5}
\mathbf{Pr}\left(\overline{A_{B,j}}\right)<e^{-\frac{\frac{n^2}{t^3}}{3\left(1-\frac{17K^2+1}{\sqrt{t}}\right) \cdot \frac{n}{t}}}
< n^{-2}.
\end{equation}

Finally, for any fixed $j\in\{1,2,\ldots,Kt\}$, let us analogously set 
$$S_j = \{v\in S: S_v\in I'_j\}.$$
Then,
\begin{eqnarray}
\mathbf{E}(|S_j|) &\leq& n \cdot \frac{\frac{1}{\sqrt{t}}\left(\frac{n}{t}+\frac{6Kn}{t\sqrt{t}}\right)}{\left(K-1\right) \left(\frac{n}{K-1} + \frac{29Kn}{\sqrt{t}}\right) \left(1-\frac{2K}{\sqrt{t}}+\frac{2K-1}{t}\right)} \nonumber\\ 
&=& \frac{1+\frac{6K}{\sqrt{t}}}{ \left(1 + \frac{29K\left(K-1\right)}{\sqrt{t}}\right) \left(1-\frac{2K}{\sqrt{t}}+\frac{2K-1}{t}\right)} \cdot \frac{n}{t\sqrt{t}} \nonumber\\ &=& 
\frac{1+\frac{6K}{\sqrt{t}}}{1+\frac{29K\left(K-1\right)-2K}{\sqrt{t}}+o\left(\frac{1}{\sqrt{t}}\right)} \cdot \frac{n}{t\sqrt{t}}\nonumber\\
 &<&  \left(1-\frac{2}{\sqrt{t}}\right) \cdot \frac{n}{t\sqrt{t}} 
 = \left(1-\frac{1}{\sqrt{t}}\right) \cdot \frac{n}{t\sqrt{t}} - \frac{n}{t^2}.\nonumber
\end{eqnarray}
Thus, by the Chernoff Bound, 
\begin{equation}\label{ProbBound6}
\mathbf{Pr}\left(\overline{A_{S,j}}\right)<e^{-\frac{\frac{n^2}{t^4}}{3\left(1-\frac{2}{\sqrt{t}}\right)\cdot\frac{n}{t\sqrt{t}}}}
< n^{-2}.
\end{equation}

By~\eqref{ProbBound1} -- \eqref{ProbBound6}
 and subadditivity of probability, all events $A_v$ -- $A'_{S,v}$ simultaneously occur with probability at least:
$$1-\left(n\cdot n^{-2}+n\cdot Kt\cdot n^{-2}+n\cdot Kt\cdot n^{-2}+n\cdot n^{-2}+n\cdot n^{-2}+n\cdot n^{-2}\right) = 1-o(1),$$
i.e., there is a choice of $x_v$, $y_v$, $M$ and $L$ (and $z_e$) fulfilling all our requirements.
\qed
\end{pf}

\subsection{Adjusting sums in $B$}

Let us fix any $x_v$, $y_v$, $M$ and $L$ consistent with the requirements of Claim~\ref{MainRandomClaim}. Let $\omega_0:E\to[K]$ be the edge weighting such that: 
$$\omega_0(e)=\left\{\begin{array}{rcl}
K &if & e\in L\\
1 &if & e\in E\smallsetminus L
\end{array}\right..$$
Let $\sigma_0(v)=Kd_L(v)+d_{E\smallsetminus L}(v)$ denote the resulting \emph{initial sum} of every vertex $v\in V$ (cf. $\sigma_v$ in Claim~~\ref{MainRandomClaim}).
These surely belong to $\bigcup_{1\leq j<Kt} I_j$.
By $A_v$, $A_{B,j}$ in Claim~~\ref{MainRandomClaim} and the definitions of $I_j$, $I'_j$,
each fixed interval $I_j$ with $j<Kt$ contains the initial sums of at most 
\begin{equation}\label{SumsInInterval}
\left(1-\frac{17K^2}{\sqrt{t}}\right)\frac{n}{t}
\end{equation}
vertices from $B$ (as these deviate by at most $ \frac{3Kn}{t\sqrt{t}}$ from their associated $S_v$). We shall move the sum of every such vertex $v$ to the next interval $I_{j+1}$ by increasing the weight of appropriate number of edges in $M$ incident with $v$ from $1$ to $2$. 
Note that by $A_{B,v}$ we can thereby increase the sum of any vertex $v$ in $B$ with $\sigma_0(v)\in I_j$ to any integer value in $I_{j+1}$.
While performing this operation we shall guarantee that these new increased sums in $I_{j+1}$ are pairwise distinct and not congruent to $0$ or $1$ modulo 
\begin{equation}\label{LambdaDefinition}
\lambda :=\left\lfloor\frac{\sqrt{t}}{8K^2}\right\rfloor.
\end{equation}
This can be achieved by a natural greedy algorithm (analysing interval after interval, starting from the last one, and increasing sums of consecutive vertices in it to first available values consistent with our requirements), as in each interval $I_{j+1}$ the number of integers which are not congruent to $0$ or $1$ modulo $\lambda$ equals at least
$$
\left(\frac{n}{t}-1-2\right) \frac{\lambda -2}{\lambda}
> \left(\frac{n}{t}-3\right) \frac{\frac{\sqrt{t}}{8K^2}-3}{\frac{\sqrt{t}}{8K^2}-1}
> \left(1-\frac{17K^2}{\sqrt{t}}\right)\frac{n}{t},
$$
that is more than in~\eqref{SumsInInterval}, i.e. than the number of vertices $v$ in $B$ with $\sigma_0(v)\in I_j$, for any interval $I_j$ with $j<Kt$. 

Denote the obtained edge weighting of $G$ by $\omega_1$, and the resulting weighted degrees by $\sigma_1(v)$ for every $v\in V$. Note that the values of $\sigma_1(v)$ are now pairwise distinct for all $v\in B$.

 \subsection{Adjusting sums in $S$}

In order to settle the \emph{final sums}, we shall design a weighting $\omega_2$ by increasing some of the weights of the edges in $S$ by at most $2$. As thus far we have only increased weights of some edges in $M$ by at most $1$, at the end of our construction, for each vertex $v\in S$ and its ultimate sum $\sigma_2(v)$ we shall have:
$$0\leq \sigma_2(v)-\sigma_0(v)\leq d_M(v)+2d_S(v)\leq \frac{2Kd(v)}{\sqrt{t}}+2\frac{d(v)}{\sqrt{t}}+2\sqrt{n} < \frac{(2K+2)n}{\sqrt{t}}+2\sqrt{n},$$ 
due to $A_{S,v}$ and $A'_{S,v}$.
This means that any vertex $v\in S$ can receive the same final sum as some $u\in S$ only if 
\begin{equation}\label{LongIntervalInS}
\sigma_0(u)\in \left(\sigma_0(v)-\frac{(2K+2)n}{\sqrt{t}}-2\sqrt{n}, \sigma_0(v)+\frac{(2K+2)n}{\sqrt{t}}+2\sqrt{n}\right).
\end{equation}
Denote the set of all such (dangerous for $v$) vertices $u$ by $D(v)$.
Note that for every $v\in V$, the number of intervals $I_j$ sufficient to cover the interval in~\eqref{LongIntervalInS} equals at most 
$$\frac{\frac{(4K+4)n}{\sqrt{t}}+4\sqrt{n}}{\frac{n}{t}}+2 < (4K+4)\sqrt{t}+3.$$
By $A_{S,j}$ and $A_v$ in Claim~\ref{MainRandomClaim} (and the definitions of $I_j$, $I'_j$), we thus obtain that for each $v\in S$,
\begin{equation}\label{D(v)Size}
\left|D(v)\right| \leq \left( (4K+4)\sqrt{t}+3\right)\left(1-\frac{1}{\sqrt{t}}\right)\frac{n}{t\sqrt{t}} 
< (4K+4)\frac{n}{t}.
\end{equation}

Before we modify the weights in $S$ we first choose for every vertex $v\in S$ a set of edges $E_v$ incident with $v$ and contained in $S$ so that $|E_v|\geq d_S(v)-1$ and so that all these sets are pairwise disjoint. We may do so by adding an auxiliary additional vertex $w$ (if necessary) to the graph $G[S]$ and joining it with an edge to every vertex with an odd degree in $G[S]$. The resulting graph $G'$ shall be Eulerian. We may thus traverse all edges of each of its components along the Eulerian trail and direct them accordingly. We then set as $E_v$ all edges in $G[S]$ outgoing from $v$. By $A'_{S,v}$ in Claim~\ref{MainRandomClaim}, for each $v\in S$ we thus have that
\begin{equation}\label{E_vSize}
\left|E_v\right|\geq \frac{d(v)}{2\sqrt{t}}-O(\sqrt{n}).
\end{equation}

We further change the weights of all edges in $S$ from $1$ to $2$ and fix any linear ordering in $S$. We then analyse one vertex after another consistently wth this ordering. While analysing a given vertex $v$ we shall be allowed to change the weights of the edges in $E_v$ in order to switch the weighted degree of $v$ to some integer $s_v$ congruent to $0$ modulo $\lambda$. Throughout the rest of the construction the sum of $v$ shall be required to be equal to this $s_v$ or $s_v+1$. 
To achieve this we admit changing the weight of any edge in $E_v$ by $1$ (if necessary) while analysing $v$. However, if $uv\in E_v$ and $u$ was analysed before $v$, we must assure that after the potential change of the weight of $uv$ the resulting sum of $u$ still belongs to $\{s_u,s_u+1\}$ (observe we shall have at least one option of modifying the weight of $uv$ by $1$ regardless of the fact if the current sum of $u$ equals $s_u$ or $s_u+1$). Note that due to such admissible modifications, we are able to attribute $v$ one of at least $|E_v|+1$ consecutive integers as its resulting sum. Since we must choose an integer congruent to $0$ modulo $\lambda$, this leaves us with at least
\begin{eqnarray}
\frac{\frac{d(v)}{2\sqrt{t}}-O(\sqrt{n})}{\lambda} &\geq& \frac{4K^2d(v)}{t}-O\left(\sqrt{\frac{n}{t}}\right)
\geq \frac{4K^2\left(\frac{n}{K-1} + \frac{29Kn}{\sqrt{t}}\right)}{t}-O\left(\sqrt{\frac{n}{t}}\right)\nonumber\\
&>& \frac{4K^2n}{t(K-1)} 
>(4K+4)\frac{n}{t} > |D(v)|\nonumber
\end{eqnarray}
options (cf.~\eqref{E_vSize}, \eqref{LambdaDefinition}, \eqref{deltaBound} and~\eqref{D(v)Size}). Out of these we may thus choose $s_v$ so that $s_v\neq s_u$ for every $u\in D(v)$ which was already analysed. We then perform admissible changes within $E_v$ so that the resulting sum of $v$ equals $s_v$ and move to the next vertex in the ordering. After analysing all vertices in $S$ we obtain our final edge weighting $\omega_2$ of $G$ and the corresponding weighted degrees $\sigma_2(v)$ for $v\in V$.

Note each edge in $S$ was first changed to $2$ within this final stage and then changed at most once by $1$ (since the sets $E_v$ are pairwise disjoint), hence $\omega_2(e)\in\{1,2,3\}$ for each such edge $e$. Since we  have not changed the weights $\omega_1(e)$ of the remaining edges, each of these is weighted $1$, $K$ or possibly $2$ (which concerns some edges in $M$).

Moreover, by the modifications performed in the last stage, every vertex $v\in S$ has its final sum $\sigma_2(v)$ distinct from $s_2(u)$ for all $u\in D(v)$. By~\eqref{LongIntervalInS}, all vertices in $S$ have thus attributed pairwise distinct sums. What is more, each of these sums is congruent to $0$ or $1$ modulo $\lambda$. As we have not modified the sums in $B$ in the last stage, these are in turn not congruent to $0$ or $1$ and remain pairwise distinct. All weighted degrees $\sigma_2(v)$, $v\in V$ are thus different, which finishes the proof of Theorem~\ref{main_theorem_constant_K}. 
\qed

\section{Concluding remarks}

Let us note that we were not striving hard to optimise the constants in the proof of Theorem~\ref{main_theorem_constant_K}. We were however trying to provide the best 
(within our approach) order of magnitude of the second term in the lower bound for $\delta$ in~\eqref{DeltaLowerInMain}.

As mentioned, due to~\eqref{LoweBoundRegular}, Theorem~\ref{main_theorem_constant_K} yields an asymptotically optimal solution of Problem~\ref{Problem29}, concerning graphs with $s(G)\leq 3$. At the same time, it also gives an asymptotically optimal solution to Problem~\ref{Problem31}, as one may easily provide a family of arbitrarily large regular graphs with $s(G)>3$ and $\delta\geq (\frac{1}{2}-o(1))n$, as indicate in particular our comments concerning~\eqref{LowerBoundAboveK} (with $K=3$).

Theorem~\ref{main_theorem_constant_K}  also yields an asymptotically optimal result related with Problem~\ref{Problem26}. In fact, we believe that small graphs are anomalies which artificially deform thresholds investigated within Problem~\ref{Problem26}, obscuring their optimal levels adequate for almost all graphs, which should be close to the ones stemming from the lower bound in~\eqref{LoweBoundRegular}. 
(In particular, in~\cite{Chartrand} it was shown that $s(G)= 5$ for the Petersen graph, while the lower bound  in~\eqref{LoweBoundRegular} yields $4$.)
Reformulating Problem~\ref{Problem26} in this vein, for every given integer $K\geq 3$ we would ask for the maximum $\alpha$ such that $s(G)\leq K$ for every graph $G$ with minimum degree $\delta \geq (1-\alpha)n$ and sufficiently large order $n$ (i.e. larger than some constant $n_K$). Then, Theorem~\ref{main_theorem_constant_K}, combined with~\eqref{LowerBoundAboveK}, yields the optimal solution of this problem for every fixed $K$, with $1-\alpha = \frac{1}{K-1}$.

Still, Theorem~\ref{main_theorem_constant_K} leaves some space for improvements, especially with respect to small graphs or more specific, accurate bounds, which seem to be quite challenging in view of the over 30 year history of these problems.

\end{document}